\documentclass{article}

\usepackage[utf8]{inputenc}
\begin{document}

\title{The Unknown Subgroup of Aut($E_8$)}

\author{Majid Butler, Sandernisha Claiborne, Tomme Denney, \\ { De'Janeke Johnson and Tianna Robinson}
\\{ }
\\{\Large McDonogh 35 High School}\\
\\{ }
{The authors were supported by the Student Research }\\
{Fellowship program at McDonogh 35 High School in New Orleans.}\\
\\{ }}

\date{We give several descriptions of the maximal subgroup 2$A_9$$<$
Aut($E_8$),
the {simplest of which is: 2A$_9$ permutes nine scale copies of $E_8$ generated by the} 
{9 $\cdot$ 240 = 2160 norm 4 vectors.}\\
\leavevmode\\
August 10, 2017}
\maketitle

\section{INTRODUCTION}

In the 150 years since it was shown to exist, the $E_8$ lattice has been one of the most studied objects in mathematics. Despite this scrutiny a significant gap in our knowledge of $E_8$ remains; the subgroup 2$A_9$, maximal in Aut($E_8$), has defied simple description. In fact the ATLAS characterizes every maximal subgroup of the simple group O$_8$$^+$(2), save one: only 2$A_9$ exhibits blank entries across the page. The purpose of our brief paper is to illuminate this one final group.

Our results can be summarized in three ways, based on the following three perspectives: looking at $E_8$ modulo 2; looking at the lattice's norm 2 vectors; and looking at the norm 4 vectors. These perspectives reveal 2$A_9$ as:

(i) the automorphisms of $E_8$ that stabilize a partition of the 135 isotropic points of $E_8$/2$E_8$ into nine disjoint isotropic 4-spaces of 15 points each;\footnote{The subgroup that fixes one of the nine 4-spaces permutes the other eight 4-spaces while it permutes the 15 points of the fixed space, so the isomorphism $A_8$ $\cong$ L$_4$(2) is visible within $E_8$.} or

(ii) the automorphisms of $E_8$ that stabilize a 9 x 15 array of coordinate frames such that

	a) each of the 120 root vectors of $E_8$ appears once in each row, as one of eight vectors in one of 15 coordinate frames, while
    
	b) each pair of orthogonal root vectors appear together in exactly one of the 135 frames; or

(iii) the automorphisms of $E_8$ that stabilize a partition of the lattice's 2160 norm 4 vectors into nine sets of 240, where each set generates a scale copy of the $E_8$ lattice itself.

\section{MODULO 2}

We begin by showing that the 135 isotropic points in an 8-space over \textbf{F}$_2$ (which admits a quadratic form of Witt defect 0) can be arranged into nine disjoint isotropic 4-spaces that are permuted by a subgroup $A_9$ $<$ O$_8$$^+$(2).

There are 270 isotropic 4-spaces in the 8-space, and O$_8$$^+$(2).2 permutes them transitively while the subgroup O$_8$$^+$(2) of index 2 permutes them in two orbits of 135 each, with stabilizer Stab($V_1$) $\cong$ $2^6$$A_8$. This stabilizer acts on the other 134 4-spaces\footnote{Unless otherwise noted, throughout the remainder of this paper the term $``$4-space$"$ will mean one of the 135 isotropic 4-spaces within the orbit in question.} in orbits of size 64 (which have trivial intersection with $V_1$) and 70 (which intersect $V_1$ in a 2-space.) If $V_2$ is a 4-space in the 64-orbit then Stab($V_1$) $\cap$ Stab($V_2$) is a subgroup isomorphic to $A_8$ that permutes the remaining 133 4-spaces in orbits of sizes 28 + 35 + 35 + 35 (which intersect \(\{\){$V_1$, $V_2$\(\}\)} in spaces of dimension \(\{\){0, 0\(\}\)}, \(\{\){0, 2\(\}\)}, \(\{\){2, 0\(\}\)} and \(\{\){2, 2\(\}\)}, respectively.)

The 28 4-spaces disjoint from both $V_1$ and $V_2$ can be labeled by the duads of eight letters. The subgroups $A_7$ $<$ $A_8$ permute the seven 4-spaces {($V_{a*}$)} = \(\{\){$V_{ab}$, $V_{ac}$, $\ldots$, $V_{ah}$\(\}\)} that share a letter, and the subgroups $A_6$ $<$ $A_8$ stabilize two 7-sets A = {($V_{a*}$)} and B = {($V_{b*}$)}, so that each $A_6$ also fixes a unique third 4-space $V_{ab}$ while permuting the six other 4-spaces in A and in B.

    We use this simple fact to construct $A_9$. We begin with two disjoint 4-spaces $V_1$ and $V_2$ (whose pointwise stabilizer is $A_8$) and we choose a subgroup $A_6$ $<$ $A_8$ that permutes six 4-spaces \(\{\){$V_{ac}$, $\ldots$, $V_{ah}$\(\}\)} and fixes a unique 4-space $V_{ab}$; note that $A_6$ extends to an $A_7$ $<$ $A_8$ that permutes the seven 4-spaces \(\{\){$V_{ab}$, $V_{ac}$, $\ldots$, $V_{ah}$\(\}\)}. But since $V_1$, $V_2$ and $V_{ab}$ are mutually disjoint, we can start instead with the $A_8$ that fixes $V_1$ and $V_{ab}$; since this group contains the above $A_6$ (which fixes $V_1$, $V_2$ and $V_{ab}$ and permutes the six 4-spaces \(\{\){$V_{ac}$, $\ldots$, $V_{ah}$\(\}\)) we can extend $A_6$ to an $A_7$ subgroup that still fixes $V_1$ and $V_{ab}$ but now permutes the seven 4-spaces \(\{\){$V_2$, $V_{ac}$, $\ldots$, $V_{ah}$\(\}\)}.\footnote{$A_6$ acts on the 28 duads of an 8-letter set in orbits of sizes 1+6+6+15 and acts on the 35 bisections in orbits of sizes 15+20, so \(\{\){$V_{ac}$, $\ldots$, $V_{ah}$\(\}\)}, an orbit of size 6, must lie in the 28-orbit for $A_8$ $\cong$ Stab($V_i$) $\cap$ Stab($V_{ab}$), i = \(\{\){1,2\(\}\)}.} Or likewise we can start with the $A_8$ that fixes $V_2$ and $V_{ab}$ to obtain an $A_7$ subgroup that permutes the seven 4-spaces \(\{\){$V_1$, $V_{ac}$, $\ldots$, $V_{ah}$\(\}\)}. These three $A_7$ subgroups together generate a subgroup of $O_8$$^+$(2) that stabilizes the nine 4-spaces \(\{\){$V_1$, $V_2$, $V_{ab}$, $V_{ac}$, $\ldots$, $V_{ah}$\(\}\)} and projects onto $A_9$; in fact this subgroup cannot be any larger than $A_9$ since otherwise its one-point stabilizer would be larger than $A_8$, which is already maximal in Stab($V_1$) $\cong$ $2^6$$A_8$.

Finally, since $V_1$ and $V_2$ are disjoint and since $A_9$ is multiply transitive on the nine 4-spaces, these nine 4-spaces are pairwise disjoint and therefore exhaust the 9$\cdot$15 = 135 isotropic vectors. 

\section{NORM 2}
                    
This arrangement of the 135 isotropic vectors into nine isotropic 4-spaces partitions the 120 root vectors of $E_8$ into nine sets of 15 coordinate frames each, with properties (ii)(a) and (ii)(b) above.

Let V be one of the nine 4-spaces and consider a 3-space W $\subset$ V. The space W$^\perp$ is 5-dimensional and contains W, so W$^\perp$/W is a 2-space with three non-zero vectors. One of these vectors completes W into V (and so is isotropic) and at least one of the other two vectors is non-isotropic, since W$^\perp$ cannot be a totally isotropic 5-space. If r is a non-isotropic vector\footnote{The third vector is isotropic and extends W to an isotropic 4-space in the other orbit of 135, leading to another description of 2A$_9$ based on the same 9 x 15 arrangement of 3-spaces.} then the eight non-isotropic vectors \(\{\){r + w, w $\in$ W\(\}\)} lift to eight different root vectors \(\{\){$r_1$, $\ldots$, $r_8$\(\}\)} of $E_8$. Since the sum $r_i$ + $r_j$ of any two of these root vectors projects to an (isotropic) vector of W, this $E_8$ vector has norm 4, so that $r_i$ and $r_j$ are orthogonal. In other words, these eight root vectors form a coordinate frame.

Furthermore, no non-isotropic vector r can lie in $W_1$$^\perp$ and $W_2$$^\perp$ for distinct 3-spaces $W_1$, $W_2$ of V since otherwise r $\in$ $W_1$$^\perp$ $\cap$ $W_2$$^\perp$ = V.  Thus the 15 3-spaces of V yield 15 frames in which each of the 120 root vectors appears only once.

The nine different 4-spaces \(\{\){$V_1$, $\ldots$, $V_9$\(\}\)} therefore give us nine different partitions of the 120 root vectors into 15 coordinate frames. No pair of root vectors \(\{\){$r_a$, $r_b$\(\}\)} may lie in two of these 135 coordinate frames (for 4-spaces $V_i$ and $V_j$, say) since then $r_a$ + $r_b$ would project to a vector in $V_i$ $\cap$ $V_j$ = \(\{\){0\(\}\)}. But since there are 120 $\cdot$ 63/2 = 3780 pairs \(\{\){$r_a$, $r_b$\(\}\)} of orthogonal root vectors, and since the 135 8-vector frames yield 135 $\cdot$ ${8 \choose 2}$  = 3780 such pairs, each pair appears in exactly one coordinate frame. Note also that since a pair \(\{\){$r_a$, $r_b$\(\}\)} of orthogonal root vectors gives rise to four norm 4 vectors \(\{\){$\pm$$r_a$$\pm$$r_b$\(\}\)}, each norm 4 vector is derived from two root vectors in 3780 $\cdot$ 4/2160 = 7 different ways.

\section{NORM 4}

Suppose we begin with the “norm 2” viewpoint: we have a 9 x 15 array of coordinate frames in which the 15 frames in each row were obtained from the 15 3-spaces of a 4-space. We use this array to partition the 2160 norm 4 vectors into nine scale copies of $E_8$, as follows. 

  Each coordinate frame of eight vectors \(\{\){$r_1$, $\ldots$, $r_8$\(\}\)} yields $\pm$56 norm 4 vectors \(\{\){$\pm$$r_i$$\pm$$r_j$\(\}\)}, and a row's 15 frames therefore yield $\pm$56$\cdot$15/7 = $\pm$120 norm 4 vectors. But since a row was obtained from an isotropic 4-space its norm 4 vectors have even inner product with each other, so we can divide the natural inner product by 2 to obtain an integral rank 8 lattice generated by $\pm$120 norm 2 vectors. If we fix one of the 15 coordinate frames in this row then the half-scale inner product turns this frame's eight root vectors \(\{\){$r_i$\(\}\)} into an orthonormal basis, in which case the $\pm$56 vectors \(\{\){$\pm$$r_i$$\pm$$r_j$\(\}\)} generate a $D_8$ lattice. The remaining $\pm$64 vectors in the other 14 frames are norm 2 elements of $D_8$$^*$, any one of which extends $D_8$ to $E_8$. And since each row contains 240 of the 2160 norm 4 vectors, the nine rows partition the norm 4 vectors into nine scale copies of $E_8$.

Finally, if the 2160 norm 4 vectors of $E_8$ are partitioned into nine double-scale copies of $E_8$, then the 240 vectors in each copy have even inner products and so form a totally isotropic subspace in $E_8$/2$E_8$. Hence these 240 vectors project onto at most 15 points in $E_8$/2$E_8$---an isotropic space cannot be larger---—and all nine copies of $E_8$ project onto no more than 135 non-zero points. So equality must hold: the nine scale copies of $E_8$ generate a partition of the 135 isotropic vectors into nine disjoint 4-spaces of 15 points each, which is exactly our original perspective.\\
\\
\Large Contacts: {majidhbutler54@gmail.com, nisha.claiborne@gmail.com,\\
tomme.denney32@gmail.com, jdejaneke@yahoo.com, and tiannarobinson71401@gmail.com}

\nocite{*}
\bibliography{main}

\begin{thebibliography}{1}

\bibitem{main}
J.H. Conway, R.T. Curtis, S.P Norton, R.A. Parker, and R.A. Wilson.
\newblock {\em An Atlas of Finite Groups}.
\newblock Oxford University Pressy, 198 Madison Avenue New York, NY 10016,
  1985.

\end{thebibliography}
\bibliographystyle{plain}

\end{document}